\documentclass{article}

\usepackage{amssymb}
\usepackage{amsthm}
\usepackage[all]{xy}

\newtheorem{thm}{Theorem}[section]
\newtheorem{prop}[thm]{Proposition}
\newtheorem{lem}[thm]{Lemma}
\newtheorem{cor}[thm]{Corollary}

\newtheorem{ithm}{Theorem}
\newtheorem{iprop}[ithm]{Proposition}

\theoremstyle{definition}
\newtheorem{dfn}[thm]{Definition}

\theoremstyle{remark}
\newtheorem{rem}{Remark}
\newtheorem*{conventions}{Conventions}
\newtheorem*{acknowledgments}{Acknowledgments}

\newcommand{\C}{\mathbb{C}}
\newcommand{\R}{\mathbb{R}}
\newcommand{\T}{\mathbb{T}}
\newcommand{\Z}{\mathbb{Z}}

\newcommand{\g}{\mathfrak{g}}

\newcommand{\F}{\mathcal{F}}
\renewcommand{\S}{\mathcal{S}}

\newcommand{\bF}{\bar{\F}}
\newcommand{\sC}{\mathcal{C}}
\newcommand{\bsC}{\bar{\mathcal{C}}}

\newcommand{\Ad}{\mathrm{Ad}}
\newcommand{\Co}{\mathrm{Co}}
\newcommand{\id}{\mathrm{id}}

\newcommand{\bCo}{\bar{\mathrm{C}}\mathrm{o}}
\newcommand{\bK}{\bar{K}}

\newcommand{\Cech}{\v{C}ech {}}
\newcommand{\Poincare}{Poincar\'e {}}
\renewcommand{\i}{\sqrt{\! - \! 1}}
\renewcommand{\d}{\partial}
\renewcommand{\epsilon}{\varepsilon}
\def\u#1{ \underline{#1} }

\def\til#1{ \tilde{#1} }

\title{Relationship between equivariant gerbes and 
gerbes over the quotient space}

\author{Kiyonori Gomi
\thanks{The author's research is supported by Research Fellowship of the Japan Society for the Promotion of Science for Young Scientists.}}

\date{}

\begin{document}

\maketitle

\begin{abstract}
By means of cohomology groups, we study relationships between equivariant gerbes with connection over a manifold with a Lie group action and gerbes with connection over the quotient space. 
\end{abstract}


\section{Introduction}

Gerbes can roughly be thought of as fiber bundles over a manifold whose fibers are categories. The notion of gerbes was originally invented by Giraud \cite{Gi} in the context of non-abelian cohomology. In \cite{Bry1}, Brylinski developed differential geometry of certain abelian gerbes (\textit{gerbes with band $\u{\T}$}, which we will simply call \textit{gerbes} from now on), and introduced the notion of connective structure and curving as connections on them. The notion of equivariant gerbes over a manifold with a Lie group action was introduced also by Brylinski \cite{Bry2,Bry1}. In the work of Sharpe \cite{Sha1,Sha2}, equivariant gerbes turned out to be useful for the geometric understanding of discrete torsions. 

Let $G$ be a Lie group acting on a smooth manifold $M$. When we can make the quotient space $M/G$ into a smooth manifold, the pull-back of a gerbe with connective structure and curving over $M/G$ by the projection $q : M \to M/G$ becomes naturally a $G$-equivariant gerbe with $G$-invariant connective structure and $G$-invariant curving over $M$. (In the sequel, we often mean a gerbe with connective structure and curving by a ``gerbe with connection'', and a $G$-equivariant gerbe with $G$-invariant connective structure and $G$-invariant curving by a ``$G$-equivariant gerbe with connection.'') The purpose of this paper is to study the relationship between $G$-equivariant gerbes with connection over $M$ and gerbes with connection over the quotient space $M/G$.

Before the study of equivariant gerbes with connection, we consider equivariant gerbes \textit{without} connection. We recall here the classification of gerbes and that of equivariant gerbes. Let $\u{\T} = \u{\T}_M$ be the sheaf of germs of smooth functions on $M$ which take its values in the unit circle $\T = \{ z \in \C |\ |z| = 1 \}$.

\begin{iprop}[Giraud \cite{Gi}] \label{iprop:classify_gerbe}
Let $M$ be a smooth manifold. The isomorphism classes of gerbes over $M$ are classified by $H^2(M, \u{\T}_M) \cong H^3(M, \Z)$.
\end{iprop}

When a Lie group $G$ acts on $M$, we denote by $G^{\bullet} \times M = \{ G^p \times M \}_{p \ge 0}$ the simplicial manifold (\cite{Se}) associated with the Lie group action. The family of sheaves $\{ \u{\T}_{G^p \times M} \}_{p \ge 0}$ gives rise to a sheaf on $G^{\bullet} \times M$, which we denote by $\u{\T}_{G^{\bullet} \times M}$. We also denote by $H^m(G^{\bullet} \times M, \u{\T}_{G^{\bullet} \times M})$ the cohomology group with coefficients $\u{\T}_{G^{\bullet} \times M}$. We note that, if $G$ is compact, then  $H^m(G^{\bullet} \times M, \u{\T}_{G^{\bullet} \times M})$ is isomorphic to the \textit{equivariant cohomology group} \cite{A-B} $H^{m+1}_G(M, \Z)$ for $m > 0$.

\begin{iprop}[Brylinski \cite{Bry2}] \label{iprop:classify_equiv_gerbe}
Let $G$ be a Lie group acting on a smooth manifold $M$. The isomorphism classes of $G$-equivariant gerbes over $M$ are classified by $H^2(G^{\bullet} \times M, \u{\T}_{G^{\bullet} \times M})$.
\end{iprop}

In the present paper, we will prove the following theorem.

\begin{ithm} \label{ithm:iso_coh}
Let $G$ be a Lie group acting on a smooth manifold $M$. We assume that the action is free and locally trivial, and that the quotient space $M/G$ is a smooth manifold in such a way that the projection map $q: M \to M/G$ is smooth. For a non-negative integer $m$ the projection map induces an isomorphism of groups
$$
q^* : H^m(M/G, \u{\T}_{M/G}) \longrightarrow 
H^m(G^{\bullet} \times M, \u{\T}_{G^{\bullet} \times M}).
$$
\end{ithm}

We can attain the assumption in Theorem \ref{ithm:iso_coh}, for example, in the case that a compact Lie group $G$ acts smoothly and freely on a finite dimensional smooth manifold $M$. However, we need not restrict ourselves to such a case only.

Combining Proposition \ref{iprop:classify_gerbe} and Proposition \ref{iprop:classify_equiv_gerbe} with Theorem \ref{ithm:iso_coh}, we directly obtain the next theorem, which is essentially known by Brylinski \cite{Bry1}.

\begin{ithm} \label{ithm:equiv_gerbe}
Let $G$ and $M$ be as in Theorem \ref{ithm:iso_coh}. 

(a) For a $G$-equivariant gerbe $\sC$ over $M$, there exists a gerbe $\bsC$ over $M/G$ whose pull-back under $q : M \to M/G$ is equivariantly isomorphic to $\sC$. 

(b) The isomorphism class of such $\bsC$ is unique.
\end{ithm}

We then consider the case of equivariant gerbes with connection by a similar method based on cohomology groups. For a non-negative integer $N$, we define a complex of sheaves $\F(N)_M$ by
$$
\u{\T}_M \stackrel{\frac{1}{2\pi\i}d\log}{\longrightarrow}
\u{A}^1_M \stackrel{d}{\longrightarrow} 
\u{A}^2_M \stackrel{d}{\longrightarrow} 
\cdots \stackrel{d}{\longrightarrow}
\u{A}^N_M \longrightarrow
0 \longrightarrow \cdots,
$$
where $\u{A}^q_M$ is the sheaf of germs of differential $q$-forms on $M$. We call the hypercohomology $H^m(M, \F(N)_M)$ the \textit{smooth Deligne cohomology} \cite{Bry1,De-F,E-V,Gaj}.

\begin{iprop}[Brylinski \cite{Bry1}] \label{iprop:classify_gerbe_conn}
Let $M$ be a smooth manifold. The isomorphism classes of gerbes with connective structure and curving over $M$ are classified by $H^2(M, \F(2)_M)$.
\end{iprop}

In \cite{Go1}, the equivariant generalization of the classification above is obtained on the basis of the work of Brylinski \cite{Bry2}. Let $G^{\bullet} \times M$ the simplicial manifold associated with the action of $G$ on $M$. The family of complexes of sheaves $\{ \F(N)_{G^p \times M} \}_{p \ge 0}$ gives rise to a complex of sheaves $\F(N)_{G^{\bullet} \times M}$ on $G^{\bullet} \times M$. 

Notice here that we have an obvious fibration $G^p \times M \to G^p \times pt$ for each $p$, where $pt$ stands for the manifold consisting of a single point. We define a subsheaf $F^1\!\u{A}^q_{G^p \times M}$ of $\u{A}^q_{G^p \times M}$ by setting $F^1\!\u{A}^q_{G^p \times M} = \pi^{-1}\u{A}^1_{G^p \times pt} \otimes \u{A}^{q-1}_{G^p \times M}$, where the tensor product is taken over $\pi^{-1}\u{A}^0_{G^p \times pt}$. Then we have a subcomplex $F^1\!\F(N)_{G^p \times M}$ of $\F(N)_{G^p \times M}$ for each $p$:
$$
\xymatrix@C=15pt@R=2pt{
\u{\T} \ar[r] &
\u{A}^1 \ar[r] &
\u{A}^2 \ar[r] &
\cdots \ar[r] &
\u{A}^{N-1} \ar[r] &
\u{A}^N \ar[r] &
0 \ar[r] &
\cdots, \\
\cup &
\cup &
\cup &
\cdots &
\cup &
\cup &
\cup & \\
0 \ar[r] &
F^1\!\u{A}^1 \ar[r] &
F^1\!\u{A}^2 \ar[r] &
\cdots \ar[r] &
F^1\!\u{A}^{N-1} \ar[r] &
F^1\!\u{A}^N \ar[r] &
0 \ar[r] &
\cdots. \\
}
$$
Thus the family $\{ F^1\!\F(N)_{G^p \times M} \}_{p \ge 0}$ gives rise to a subcomplex $F^1\!\F(N)_{G^{\bullet} \times M}$ of $\F(N)_{G^{\bullet} \times M}$. We define a complex of sheaves $\bF(N)_{G^{\bullet} \times M}$ by taking the quotient: $\bF(N)_{G^{\bullet} \times M} = \F(N)_{G^{\bullet} \times M} / F^1\!\F(N)_{G^{\bullet} \times M}$. In \cite{Go1}, the hypercohomology groups $H^m(G^{\bullet} \times M, \bF(N)_{G^{\bullet} \times M})$ are called ``\textit{equivariant smooth Deligne cohomology groups}.'' In the sequel, we will omit the subscripts of $F^1\!\F(N)_{G^{\bullet} \times M}$, $\F(N)_{G^{\bullet} \times M}$ and $\bF(N)_{G^{\bullet} \times M}$.

\begin{iprop}[\cite{Go1}] \label{iprop:classify_equiv_gerbe_conn}
Let $G$ be a Lie group acting on a smooth manifold $M$. The isomorphism classes of $G$-equivariant gerbes with $G$-invariant connective structure and $G$-invariant curving over $M$ are classified by $H^2(G^{\bullet} \times M, \bF(2))$.
\end{iprop}

We will prove the following key theorem in this paper.

\begin{ithm} \label{ithm:iso_Deligne_coh}
Let $G$ and $M$ be as in Theorem \ref{ithm:iso_coh}. For a non-negative integer $N$ the projection map $q : M \to M/G$ induces an isomorphism of groups
$$
q^* : H^N(M/G, \F(N)) \longrightarrow 
H^N(G^{\bullet} \times M, \F(N)).
$$
\end{ithm}

To relate $H^m(G^{\bullet} \times M, \F(N))$ and $H^m(G^{\bullet} \times M, \bF(N))$, we consider the long exact sequence induced by the following short exact sequence of complex of sheaves on the simplicial manifold $G^{\bullet} \times M$:
$$
0 \longrightarrow
F^1\!\F(N) \longrightarrow
\F(N) \longrightarrow
\bF(N) \longrightarrow 0.
$$
We denote the Bockstein homomorphism in the long exact sequence by
$$
\beta : H^m(G^{\bullet} \times M, \bF(N)) \longrightarrow 
H^{m+1}(G^{\bullet} \times M, F^1\!\F(N)).
$$

Now Proposition \ref{iprop:classify_gerbe_conn} and Proposition \ref{iprop:classify_equiv_gerbe_conn} lead to the following result on the relationship between $G$-equivariant gerbes with connection over $M$ and gerbes with connection over the quotient space $M/G$.

\begin{ithm} \label{ithm:equiv_gerbe_conn}
Let $G$ and $M$ be as in Theorem \ref{ithm:iso_coh}, $(\sC, \Co, K)$ a $G$-equivariant gerbe with $G$-invariant connective structure and $G$-invariant curving over $M$, and $c \in H^2(G^{\bullet} \times M, \bF(2))$ the class corresponding to the equivariant isomorphism class of $(\sC, \Co, K)$. 

(a) There exists a gerbe with connective structure and curving $(\bsC, \bCo, \bK)$ over $M/G$ whose pull-back under the projection $q : M \to M/G$ is equivariantly isomorphic to $(\sC, \Co, K)$ if and only if $\beta(c) = 0$ in $H^3(G^{\bullet} \times M, F^1\!\F(2))$.

(b) The isomorphism classes of such $(\bsC, \bCo, \bK)$ are in one to one correspondence with $Coker\{ \beta : H^1(G^{\bullet} \times M, \bF(2)) \to H^2(G^{\bullet} \times M, F^1\!\F(2)) \}$.
\end{ithm}

By expressing $H^2(G^{\bullet} \times M, F^1\!\F(2))$ in a more accessible form, we can obtain a condition for the isomorphism class of such $(\bsC, \bCo, \bK)$ as in Theorem \ref{ithm:equiv_gerbe_conn} (a) to be unique. To the contrary, by a simple computation, we can find a case in which there indeed exist distinct isomorphism classes of such $(\bsC, \bCo, \bK)$.

\bigskip

The outline of the present paper is as follows. In Section \ref{sec:SDC}, we briefly recall the smooth Deligne cohomology groups. In Section \ref{sec:EDC} we review the definition of the equivariant smooth Deligne cohomology group $H^m(G^{\bullet} \times M, \bF(N))$ and some facts on it. In Section \ref{sec:key_theorem}, we prove Theorem \ref{ithm:iso_coh} (Theorem \ref{thm:quotient_simplicial_uT}) and Theorem \ref{ithm:iso_Deligne_coh} (Theorem \ref{thm:Deligne_coh_quotient}). Section \ref{sec:ecbc} deals with relationships between equivariant circle bundles and circle bundles over the quotient space. Though the relationships can be studied directly, we use the results in Section \ref{sec:key_theorem} which have an advantage of generalization. In Section \ref{sec:egc}, we study the main subject of this paper: the relationship between equivariant gerbes with connection and gerbes with connection over the quotient space.

To save pages, we mainly follow the terminologies in \cite{Bry1,Bry2,Bry-M}, and drop the definition of (equivariant) gerbes.

\bigskip

\begin{conventions}
Throughout this paper, we make a convention that a ``smooth manifold'' means a paracompact smooth manifold modeled on a topological vector space which is Hausdorff and locally convex. We also assume the existence of a partition of unity. Examples of such a manifold cover not only all the finite dimensional smooth manifolds, but also a sort of infinite dimensional manifolds. (An example of the infinite dimensional case is the loop space of a finite dimensional smooth manifold. See \cite{Bry1} for detail.) 

We also make a convention that a ``Lie group'' means a Lie group whose underlying smooth manifold is of the type above. When a Lie group $G$ acts on a smooth manifold $M$, we denote the action by juxtaposition: we write $gx \in M$ for $g \in G$ and $x \in M$. The unit element of $G$ is denoted by $e \in G$.

We usually work in the smooth category, so functions, differential forms, etc. are assumed to be smooth.
\end{conventions}


\section{Smooth Deligne cohomology}
\label{sec:SDC}

We here recall ordinary smooth Deligne cohomology groups \cite{Bry1,De-F,E-V,Gaj}. 

\subsection{Smooth Deligne cohomology groups}

Let $M$ be a smooth manifold. We denote by $\u{\T}_M$ the sheaf of germs of functions with values in the unit circle $\T = \{ u \in \C |\ |u| = 1 \}$. For a non-negative integer $q$, we denote by $\u{A}^q_M$ the sheaf of germs of $\R$-valued differential $q$-forms on $M$. 

\begin{dfn}
Let $N$ be a non-negative integer. 

(a) We define the \textit{smooth Deligne complex} $\F(N)_M$ to be the following complex of sheaves on $M$:
$$
\F(N)_M : \
\u{\T}_M \stackrel{\frac{1}{2\pi\i}d\log}{\longrightarrow}
\u{A}^1_M \stackrel{d}{\longrightarrow} 
\u{A}^2_M \stackrel{d}{\longrightarrow} 
\cdots \stackrel{d}{\longrightarrow}
\u{A}^N_M \longrightarrow
0 \longrightarrow \cdots,
$$
where $\u{\T}_M$ is located at degree 0 in the complex. 

(b) The \textit{smooth Deligne cohomology group} $H^p(M, \F(N)_M)$ is defined to be the hypercohomology group of the smooth Deligne complex.
\end{dfn}

We often omit the subscripts of $\u{\T}_M$, $\u{A}^q_M$ and $\F(N)_M$.

\begin{rem}
Let $\Z(N)_D^\infty$ be a complex of sheaves given by
$$
\Z(N)_D^\infty : \
\Z \stackrel{i}{\longrightarrow}
\u{A}^0 \stackrel{d}{\longrightarrow} 
\u{A}^1 \stackrel{d}{\longrightarrow} 
\u{A}^2 \stackrel{d}{\longrightarrow}
\cdots \stackrel{d}{\longrightarrow}
\u{A}^{N-1} \longrightarrow
0 \longrightarrow \cdots,
$$
where $\Z$ is regarded as a constant sheaf. The smooth Deligne cohomology often refers to the hypercohomology $H^p(M, \Z(N)_D^\infty)$. Since $\Z(N)_D^\infty$ is quasi-isomorphic to $\F(N-1)$ under a shift of degree, we have $H^p(M, \Z(N)_D^\infty) \cong H^{p-1}(M, \F(N-1))$.
\end{rem}

Recall the following classification of principal $\T$-bundles and gerbes. (In this paper, a ``gerbe'' means a \textit{gerbe with band $\u{\T}$} \cite{Bry1,Bry-M,Gi}.)

\begin{prop} \label{prop:classify}
Let $M$ be a smooth manifold.

(a)(Kostant \cite{Ko}, Weil \cite{We}) The isomorphism classes of principal $\T$-bundles (Hermitian line bundles) over $M$ are classified by $H^1(M, \u{\T}) \cong H^2(M, \Z)$.

(b)(Giraud \cite{Gi}) The isomorphism classes of gerbes over $M$ are classified by $H^2(M, \u{\T}) \cong H^3(M, \Z)$.
\end{prop}

By using the smooth Deligne cohomology groups, we obtain the following generalization of the proposition above.

\begin{prop}[Brylinski \cite{Bry1}] \label{prop:classify_connection}
Let $M$ be a smooth manifold.

(a) The isomorphism classes of principal $\T$-bundles with connection over $M$ are classified by $H^1(M, \F(1))$.

(b) The isomorphism classes of gerbes with connective structure and curving over $M$ are classified by $H^2(M, \F(2))$.
\end{prop}

We omit the proofs of Proposition \ref{prop:classify} and Proposition \ref{prop:classify_connection}, and refer the reader to \cite{Bry1}.


\section{Equivariant smooth Deligne cohomology}
\label{sec:EDC}

This section is a short summary of \cite{Go1}. We introduce equivariant smooth Deligne cohomology groups, and state the classification of equivariant circle bundles (with connection) and equivariant gerbes (with connection).

\subsection{Simplicial manifolds associated to group actions}

Let $G$ be a Lie group acting on a smooth manifold $M$ by left. Then we have a simplicial manifold $G^{\bullet} \times M = \{ G^p \times M \}_{p \ge 0}$, where the face maps $\d_i : G^{p+1} \times M \to G^{p} \times M$, $(i = 0, \ldots p+1)$ are given by
$$
\d_i(g_1, \ldots, g_{p+1}, x) =
\left\{
\begin{array}{ll}
(g_2, \ldots, g_{p+1}, x), & i = 0 \\
(g_1, \ldots, g_{i-1}, g_i g_{i+1}, g_{i+2}, \ldots, g_{p+1}, x), & i = 1, \ldots, p \\
(g_1, \ldots, g_p, g_{p+1} x), & i = p + 1,
\end{array}
\right.
$$
and the degeneracy maps $s_i : G^p \times M \to G^{p+1} \times M$, $(i = 0, \ldots p)$ by
$$
s_i(g_1, \ldots, g_p, x) = (g_1, \ldots, g_i, e, g_{i+1}, \ldots, g_p, x).
$$
These maps obey the following relations:
\begin{eqnarray}
\d_i \circ \d_j & = & \d_{j-1} \circ \d_i, \quad (i < j), \label{rel1} \\
s_i \circ s_j & = & s_{j+1} \circ s_i, \quad (i \le j), \label{rel2} \\
\d_i \circ s_j & = &
\left\{
\begin{array}{ll}
s_{j-1} \circ \d_i, & \quad (i < j), \\
\id, & \quad (i = j, j+1), \\
s_j \circ \d_{i-1}, &\quad (i > j+1).
\end{array}
\right. \label{rel3}
\end{eqnarray}

To a simplicial manifold, we can associate a topological space called the \textit{realization} \cite{Gaj,Se}. The realization of $G^{\bullet} \times M$ is identified with the homotopy quotient (\cite{A-B}): $|G^{\bullet} \times M| \cong (EG \times M) /G$, where $EG$ is the total space of the universal bundle for $G$. This can be seen by the fact that $EG$ is obtained as the realization of $G^{\bullet} \times G$, where $G$ acts on itself by the left translation. 

Note that the classifying space $BG$ is also obtained as the realization of $G^{\bullet} \times pt$, where $pt$ is the space consisting of a single point on which $G$ acts trivially. We denote by $\pi : G^{\bullet} \times M \to G^{\bullet} \times pt$ the map of simplicial manifolds given by the projection $\pi : G^p \times M \to G^p \times pt$.

\subsection{Equivariant smooth Deligne cohomology}

We here explain briefly the notion of a \textit{sheaf on a simplicial manifold} (a \textit{simplicial sheaf}, for short) \cite{De}. Let $G^{\bullet} \times M$ be the simplicial manifold associated to an action of a Lie group $G$ on a smooth manifold $M$. We define a \textit{simplicial sheaf} on $G^{\bullet} \times M$ to be a family of sheaves $\S^{\bullet} = \{ \S^p \}_{p \ge 0}$, where $\S^p$ is a sheaf on $G^p \times M$ such that homomorphisms $\til{\d}_i : \d_i^{-1} \S^p \to \S^{p+1}$ and $\til{s}_i : s_i^{-1} \S^{p+1} \to \S^p$ obeying the same relations as (\ref{rel1}), (\ref{rel2}) and (\ref{rel3}) are specified. 

For example, the family $\{ \u{\T}_{G^p \times M} \}_{p \ge 0}$ gives rise to a simplicial sheaf on $G^{\bullet} \times M$, which we denote by $\u{\T}_{G^{\bullet} \times M}$ or $\u{\T}$.

\smallskip

Let $\S^{\bullet} = \{ \S^p\}_{p \ge 0}$ be a simplicial sheaf on $G^\bullet \times M$. For each $p$, let $(I^{p, *}, \delta)$ be an injective resolution of $\S^{p}$. We call $I^{*, *}$ an injective resolution of $\S^{\bullet}$. The homomorphism $\til{\d}_i : \d_i^{-1} \S^p \to \S^{p+1}$ induces a homomorphism $\d_i^* : \Gamma(G^p \times M, I^{p,q}) \to \Gamma(G^{p+1} \times M, I^{p+1,q})$. Combining these homomorphisms, we define a homomorphism $\d : \Gamma(G^p \times M, I^{p,q}) \to \Gamma(G^{p+1} \times M, I^{p+1,q})$ to be $\d = \sum_{j=0}^{p+1} (-1)^j\d_j^*$. This homomorphism satisfies $\d \circ \d = 0$, because of (\ref{rel1}). We define $H^*(G^{\bullet} \times M, \S^{\bullet})$, the cohomology with coefficients the simplicial sheaf $\S^{\bullet}$, to be the total cohomology of the double complex $(\Gamma(G^i \times M, I^{i, j}), \d, \delta)$.

\smallskip

The notion of a complex of simplicial sheaves and of its hypercohomology are defined in a similar fashion.

\begin{dfn}
Let $G$ be a Lie group acting on a smooth manifold $M$. We define a complex of simplicial sheaves $\F(N)_{G^{\bullet} \times M}$ on $G^{\bullet} \times M$ by the family of smooth Deligne complex $\{ \F(N)_{G^p \times M} \}_{p \ge 0}$, where the homomorphisms $\til{\d}_i : \d_i^{-1} \F(N)_{G^p \times M} \to \F(N)_{G^{p+1} \times M}$ and $\til{s}_i : s_i^{-1} \F(N)_{G^{p+1} \times M} \to \F(N)_{G^p \times M}$ are the natural ones.
\end{dfn}

Here we consider the smooth Deligne complex $\F(N)_{G^i \times M}$ on $G^i \times M$ for an $i$ fixed. We have an obvious fibration $\pi : G^i \times M \to G^i \times pt$. For a positive integer $p$, we define a subsheaf $F^p\!\u{A}^q_{G^i \times M}$ of $\u{A}^q_{G^i \times M}$ by setting $F^p\!\u{A}^q_{G^i \times M} = \pi^{-1}\u{A}^p_{G^i \times pt} \otimes \u{A}^{q-p}_{G^i \times M}$, where the tensor product is taken over $\pi^{-1}\u{A}^0_{G^i \times pt}$.

For an open subset $U \subset G^i \times M$, the group $F^p\!\u{A}^q_{G^i \times M}(U)$ consists of those $q$-forms $\omega$ on $U$ satisfying $\iota_{V_1} \cdots \iota_{V_{q-p+1}} \omega = 0$ for tangent vectors $V_1, \ldots, V_{q-p+1}$ at $x \in U$ such that $\pi_*V_k = 0$. If $\{ g_j \}$ and $\{ x_k \}$ are systems of local coordinates of $G$ and $M$ respectively, then the $q$-form $\omega$ has a local expression
$$
\omega = 
\sum_{r \ge p} 
\sum_{{j_1, \ldots, j_r}, \atop {k_1, \ldots, k_{q-r}}}
f_{J, K}(g, x) 
dg_{j_1} \wedge \cdots \wedge dg_{j_r} \wedge
dx_{k_1} \wedge \cdots \wedge dx_{k_{q-r}}.
$$

Since we have a filtration $\u{A}^q_{G^i \times M} \supset F^1\!\u{A}^q_{G^i \times M} \supset \cdots \supset F^q\!\u{A}^q_{G^i \times M} \supset 0$, the smooth Deligne complex $\F(N)_{G^i \times M}$ admits the following filtration:
$$
\xymatrix@C=15pt@R=2pt{
\u{\T} \ar[r] &
\u{A}^1 \ar[r] &
\u{A}^2 \ar[r] &
\cdots \ar[r] &
\u{A}^{N-1} \ar[r] &
\u{A}^N \ar[r] &
0 \ar[r] &
\cdots, \\
\cup &
\cup &
\cup &
\cdots &
\cup &
\cup &
\cup & \\
0 \ar[r] &
F^1\!\u{A}^1 \ar[r] &
F^1\!\u{A}^2 \ar[r] &
\cdots \ar[r] &
F^1\!\u{A}^{N-1} \ar[r] &
F^1\!\u{A}^N \ar[r] &
0 \ar[r] &
\cdots, \\
\cup &
\cup &
\cup &
\cdots &
\cup &
\cup &
\cup & \\
0 \ar[r] &
0 \ar[r] &
F^2\!\u{A}^2 \ar[r] &
\cdots \ar[r] &
F^2\!\u{A}^{N-1} \ar[r] &
F^2\!\u{A}^N \ar[r] &
0 \ar[r] &
\cdots, \\
\cup &
\cup &
\cup &
\cdots &
\cup &
\cup &
\cup & \\
\vdots &
\vdots &
\vdots &
\cdots &
\vdots &
\vdots &
\vdots &
\cdots, \\
\cup &
\cup &
\cup &
\cdots &
\cup &
\cup &
\cup & \\
0 \ar[r] &
0 \ar[r] &
0 \ar[r] &
\cdots \ar[r] &
0 \ar[r] &
F^N\!\u{A}^N \ar[r] &
0 \ar[r] &
\cdots, \\
\cup &
\cup &
\cup &
\cdots &
\cup &
\cup &
\cup & \\
0 \ar[r] &
0 \ar[r] &
0 \ar[r] &
\cdots \ar[r] &
0 \ar[r] &
0 \ar[r] &
0 \ar[r] &
\cdots.
}
$$
We denote this filtration by 
$$
\F(N)_{G^i \times M} \supset
F^1\!\F(N)_{G^i \times M} \supset
F^2\!\F(N)_{G^i \times M} \supset \cdots \supset
F^N\!\F(N)_{G^i \times M} \supset 0.
$$

For each $i$, we also define a complex of sheaves $\bF(N)_{G^i \times M}$ by taking the quotient: $\bF(N)_{G^i \times M} = \F(N)_{G^i \times M} / F^1\!\F(N)_{G^i \times M}$. If we introduce the sheaf of germs of \textit{relative} $q$-forms with respect to $G^i \times M \to G^i \times pt$ by $\u{A}^q_{rel} = \u{A}^q_{G^i \times M} / F^1\!\u{A}^q_{G^i \times M}$, then the complex $\bF(N)_{G^i \times M}$ is expressed as 
$$
\bF(N) : \
\u{\T} \stackrel{\frac{1}{2\pi\i}d\log}{\longrightarrow}
\u{A}^1_{rel} \stackrel{d}{\longrightarrow} 
\u{A}^2_{rel} \stackrel{d}{\longrightarrow} 
\cdots \stackrel{d}{\longrightarrow}
\u{A}^N_{rel} \longrightarrow
0 \longrightarrow \cdots.
$$

\begin{dfn}
Let $G$ be a Lie group acting on a smooth manifold $M$.

(a) We define a subcomplex $F^1\!\F(N)_{G^{\bullet} \times M}$ of $\F(N)_{G^{\bullet} \times M}$ by the family $\{ F^1\!\F(N)_{G^i \times M} \}_{i \ge 0}$.

(b) We define a complex of simplicial sheaves $\bF(N)_{G^{\bullet} \times M}$ on $G^{\bullet} \times M$ by the family $\{ \bF(N)_{G^i \times M} \}_{i \ge 0}$.
\end{dfn}

We can also define $\bF(N)_{G^{\bullet} \times M}$ by the quotient $\F(N)_{G^{\bullet} \times M} / F^1\!\F(N)_{G^{\bullet} \times M}$. As is clear, if the topology of $G$ is discrete, then $F^1\!\F(N)_{G^{\bullet} \times M} = 0$, so that we have $\F(N)_{G^{\bullet} \times M} = \bF(N)_{G^{\bullet} \times M}$.

\begin{dfn}[\cite{Go1}] \label{dfn:EDC}
Let $G$ be a Lie group acting on a smooth manifold $M$. We define the \textit{$G$-equivariant smooth Deligne cohomology group} of $M$ to be the hypercohomology group $H^m(G^{\bullet} \times M, \bF(N)_{G^{\bullet} \times M})$ of the complex of simplicial sheaves $\bF(N)_{G^{\bullet} \times M}$ on $G^{\bullet} \times M$.
\end{dfn}

From now on, we omit the subscripts of $\F(N)_{G^i \times M}$, $\F(N)_{G^{\bullet} \times M}$, etc. 

\begin{rem}
When $G$ is a finite group, the hypercohomology $H^m(G^{\bullet} \times M, \F(N))$ is introduced in the work of Lupercio and Uribe \cite{L-U} as the \textit{Deligne cohomology group for the orbifold $M/G$}. Since the topology of $G$ is discrete, the cohomology $H^m(G^{\bullet} \times M, \bF(N))$ coincides with $H^m(G^{\bullet} \times M, \F(N))$.
\end{rem}

By definition, the hypercohomology $H^m(G^{\bullet} \times M, \bF(N))$ is given in the following way. Let $I^{*, *, *}$ be an injective resolution of $\bF(N)$, that is, $I^{i, *, *}$ is an injective resolution of the complex of sheaves $\bF(N)$ on $G^i \times M$:
$$
\xymatrix{
\vdots &
\vdots &
&
\vdots &
\vdots \\
I^{i, 1, 0} \ar[r]^{\til{d}} \ar[u]^{\delta} &
I^{i, 1, 1} \ar[r]^{\til{d}} \ar[u]^{\delta} &
\cdots \ar[r]^-{\til{d}} &
I^{i, 1, N} \ar[r]^{\til{d}} \ar[u]^{\delta} &
I^{i, 1, N+1} \ar[r]^{\til{d}} \ar[u]^{\delta} &
\cdots \\
I^{i, 0, 0} \ar[r]^{\til{d}} \ar[u]^{\delta} &
I^{i, 0, 1} \ar[r]^{\til{d}} \ar[u]^{\delta} &
\cdots \ar[r]^-{\til{d}} &
I^{i, 0, N} \ar[r]^{\til{d}} \ar[u]^{\delta} &
I^{i, 0, N+1} \ar[r]^{\til{d}} \ar[u]^{\delta} &
\cdots \\
\u{\T} \ar[r]^{\til{d}} \ar[u] &
\u{A}^1_{rel} \ar[r]^{\til{d}} \ar[u] &
\cdots \ar[r]^-{\til{d}} &
\u{A}^N_{rel} \ar[r]^{\til{d}} \ar[u] &
0 \ar[r]^{\til{d}} \ar[u] &
\cdots.
}
$$
We put $K^{i, j, k} = \Gamma(G^i \times M, I^{i, j, k})$. The injective resolution induces coboundary operators $\delta : K^{i, j, k} \to K^{i, j+1, k}$ and $\til{d} : K^{i, j, k} \to K^{i, j, k+1}$. The homomorphism $\til{\d}_l : \d_l^{-1}\bF(N)_{G^i \times M} \to \bF(N)_{G^{i+1} \times M}$ induces $\d_l^* : K^{i, j, k} \to K^{i+1, j, k}$. If we define $\d : K^{i, j, k} \to K^{i+1, j, k}$ by $\d = \sum_{l=0}^{i+1} (-1)^l\d_l^*$, then we have $\d \circ \d = 0$ by the relation (\ref{rel1}). Since $\d$ commutes with both $\delta$ and $\til{d}$, we have a triple complex $(K^{i, j, k}, \d, \delta, \til{d})$. For $\oplus_{m = i+j+k} K^{i, j, k}$, the total coboundary operator is defined by $D = \d + (-1)^i \delta + (-1)^{i+j} \til{d}$ on the component $K^{i, j, k}$. The cohomology of this total complex is $H^m(G^{\bullet} \times M, \bF(N))$.

\medskip

Let $\{ F^pK \}_{p = 0, 1, \ldots}$ be a filtration of the triple complex $K^{*, *, *}$ given by $F^{p}K = \oplus_{i \ge p} K^{i, *, *}$. This provides us a spectral sequence converging to the graded quotient of $H^m(G^{\bullet} \times M, \bF(N))$ with respect to the filtration. The $E_1$-terms are 
\begin{eqnarray}
E^{p, q}_1 = H^q(G^p \times M, \bF(N)),
\label{E_1:Type_Ia}
\end{eqnarray}
and the differential $d_1 : E_1^{p, q} \to E_1^{p+1, q}$ is $\d = \sum_{l=0}^{p+1} (-1)^l \d_l^*$. Note that $E^{0, q}_1$ coincides with the ordinary smooth Deligne cohomology $H^q(M, \F(N))$. 

\begin{lem} \label{lem:G_trivial}
If $G = \{ e \}$, then $H^m(G^{\bullet} \times M, \bF(N)) \cong H^m(M, \F(N))$.
\end{lem}

\begin{proof}
We use the spectral sequence (\ref{E_1:Type_Ia}). The natural identification $\{ e \}^p \times M = M$ implies that $E^{p, q}_1 = E^{0, q}_1$ for all $p$ and $q$. Under this identification, $d_1 = 0$ if $p$ is even, and $d_1 = \id$ if $p$ is odd. Thus, the spectral sequence degenerates at $E_2$, and we obtain $H^q( \{ e \}^{\bullet} \times M, \bF(N)) = E_2^{0, q} = H^q(M, \F(N))$.
\end{proof}

\subsection{The classification of equivariant circle bundles and equivariant gerbes}

In \cite{Bry2}, Brylinski classified equivariant principal $\T$-bundles and equivariant gerbes by means of the cohomology $H^m(G^{\bullet} \times M, \u{\T}) \cong H^m(G^{\bullet} \times M, \bF(0))$.

\begin{prop}[Brylinski \cite{Bry2}] \label{prop:classify_equiv}
Let $G$ be a Lie group acting on a smooth manifold $M$. 

(a) The isomorphism classes of $G$-equivariant principal $\T$-bundles over $M$ are classified by $H^1(G^{\bullet} \times M, \u{\T})$.

(b) The isomorphism classes of $G$-equivariant gerbes over $M$ are classified by $H^2(G^{\bullet} \times M, \u{\T})$.
\end{prop}

\begin{rem} \label{rem_EC}
Let $EG$ be the total space of the universal $G$-bundle. For a smooth manifold $M$ with a $G$-action, the \textit{equivariant cohomology group} \cite{A-B} is often defined by $H_G^m(M, \Z) = H^m((EG \times M)/G, \Z)$, where $G$ acts on $EG \times M$ diagonally. If $G$ is compact and $m$ is a positive integer, then $H^m(G^{\bullet} \times M, \u{\T})$ is isomorphic to $H_G^{m+1}(M, \Z)$. (See \cite{Bry2}.)
\end{rem}

The equivariant smooth Deligne cohomology $H^m(G^{\bullet} \times M, \bF(N))$ allows one to have the following generalization of the proposition above.

\begin{prop}[\cite{Go1}] \label{prop:classify_equiv_connection}
Let $G$ be a Lie group acting on a smooth manifold $M$. 

(a) The isomorphism classes of $G$-equivariant principal $\T$-bundles with $G$-invariant connection over $M$ are classified by $H^1(G^{\bullet} \times M, \bF(1))$.

(b) The isomorphism classes of $G$-equivariant gerbes with $G$-invariant connective structure and $G$-invariant curving over $M$ are classified by $H^2(G^{\bullet} \times M, \bF(2))$.
\end{prop}

Proposition \ref{prop:classify_equiv} and Proposition \ref{prop:classify_equiv_connection} are shown by using a \Cech cohomology description of $H^m(G^{\bullet} \times M, \u{\T})$ and $H^m(G^{\bullet} \times M, \bF(N))$. See \cite{Bry2,Go1} for detail.

\section{Key theorem}
\label{sec:key_theorem}

\subsection{The cohomology on the quotient space}

Let $G$ be a Lie group acting on a smooth manifold $M$. We endow the quotient space $M/G$ with the quotient topology, so that the natural projection map $q : M \to M/G$ is a continuous map.

\begin{lem} \label{lem:Leray_type_ss_uT}
Let $G$ be a Lie group acting on a smooth manifold $M$. There exists a spectral sequence converging to a graded quotient of the cohomology group $H^{p+q}(G^{\bullet} \times M, \u{\T})$ with its $E_2$-term given by
$$
E^{p, q}_2 = H^p(M/G, \u{\mathcal{X}}^q),
$$
where $\u{\mathcal{X}}^q$ is the sheaf on $M/G$ associated with the presheaf given by the assignment of $H^q(G^{\bullet} \times q^{-1}(V), \u{\T})$ to an open subset $V \subset M/G$.
\end{lem}

\begin{proof}
The spectral sequence is a sort of the Leray spectral sequence \cite{Bry1}. Let $q : G^{\bullet} \times M \to \{ e \}^{\bullet} \times (M/G)$ be the simplicial map induced by the projection $q: M \to M/G$, and $I^{*, *}$ an injective resolution of the simplicial sheaf $\u{\T}$ on $G^{\bullet} \times M$. We denote by $q_*I^{i, j}$ the direct image of $I^{i, j}$ under the projection $q : G^{i} \times M \to \{ e \}^i \times (M/G)$. Since we can identify $\{ e \}^{i} \times (M/G)$ with $M/G$, we have a double complex of sheaves $q_*I^{*, *}$ on $M/G$. We compute the hypercohomology of the complex of sheaves $q_*I^* = \oplus_{*=i+j}q_*I^{i, j}$ on $M/G$ in two ways. We take an injective resolution $J^{*, *}$ of the complex of sheaves $q_*I^*$:
$$
\xymatrix{
\vdots &
\vdots &
\vdots \\
J^{1, 0} \ar[r] \ar[u] &
J^{1, 1} \ar[r] \ar[u] &
J^{1, 2} \ar[r] \ar[u] &
\cdots \\
J^{0, 0} \ar[r] \ar[u] &
J^{0, 1} \ar[r] \ar[u] &
I^{0, 2} \ar[r] \ar[u] &
\cdots \\
q_*I^0 \ar[r] \ar[u] &
q_*I^1 \ar[r] \ar[u] &
q_*I^2 \ar[r] \ar[u] &
\cdots.
}
$$

On the one hand, the filtration ${}'\!F^p = \oplus_{j \ge p}J^{*, j}$ induces a spectral sequence converging to a graded quotient of $H^{p+q}(M/G, q_*I^*)$ whose $E_1$-terms are
$$
{}'\!E_1^{p, q} = 
\left\{
\begin{array}{cl}
\Gamma(M/G, q_*I^p), & \quad (q = 0), \\
0, & \quad (q > 0).
\end{array}
\right.
$$
Since $\Gamma(M/G, q_*I^p) = \Gamma(M/G, \oplus_{p=i+j}q_*I^{i, j}) = \oplus_{p = i+j}\Gamma(G^i \times M, I^{i, j})$, the $E_2$-terms become
$$
{}'\!E_2^{p, q} = 
\left\{
\begin{array}{cl}
H^p(G^{\bullet} \times M, \u{\T}), & \quad (q = 0), \\
0, & \quad (q > 0).
\end{array}
\right.
$$
Thus, the spectral sequence degenerates at $E_2$, and yields an isomorphism $H^m(M/G, q_*I^*) \cong H^m(G^{\bullet} \times M, \u{\T})$.

On the other hand, the filtration ${}''\!F^p = \oplus_{i \ge p}J^{i, *}$ gives another spectral sequence converging to a graded quotient of $H^m(M/G, q_*I^*)$. Its $E_2$-terms are given by ${}''\!E^{p, q}_2 = H^p(M/G, \u{H}^q(q_*I^*))$, where $\u{H}^q(q_*I^*)$ is the $q$th cohomology sheaf of $q_*I^*$, namely, the sheaf associated with the presheaf $V \mapsto H^q(V, q_*I^*) = H^q(G^{\bullet} \times q^{-1}(V), \u{\T})$. Since $\u{\mathcal{X}}^q = \u{H}^q(q_*I^*)$ by definition, we obtain the spectral sequence in this lemma.
\end{proof}

\begin{lem} \label{lem:cohomology_vanishing_uT}
Let $V$ be a contractible smooth manifold. If $G$ acts on $G \times V$ by the left translation on $G$ and by the trivial action on $V$, then we have
$$
H^m(G^{\bullet} \times (G \times V), \u{\T}) =
\left\{
\begin{array}{cl}
C^{\infty}(V, \T), & \quad (m = 0), \\
0, & \quad (m > 0),
\end{array}
\right.
$$
where $C^{\infty}(V, \T)$ is the group of smooth $\T$-valued functions on $V$.
\end{lem}

\begin{proof}
We use the spectral sequence (\ref{E_1:Type_Ia}):
\begin{eqnarray*}
E_1^{p, q} & = & H^q(G^p \times (G \times V), \u{\T}), \\
E_2^{p, q} & = & H^p( H^q(G^* \times (G \times V), \u{\T}), \d).
\end{eqnarray*}
We define a map $\phi_p : G^p \times (G \times V) \to G^{p+1} \times (G \times V)$ by 
$$
\phi_p (g_1, \ldots, g_p, h, x) = (g_1, \ldots, g_p, h, e, x).
$$
If $p > 0$, then they obey $\d_i \circ \phi_p = \phi_{p-1} \circ \d_i$ for $i < p+1$ and $\d_{p+1} \circ \phi_p = \id$. Thus, if $c \in E_1^{p, q}$ is a class such that $\d c = 0$, then $\phi_{p-1}^*c \in E_1^{p-1, q}$ satisfies $\d(\phi_{p-1}^*c) = (-1)^pc$. Hence $E^{p, q}_2 = 0$ for all $p > 0$ and $q$. The spectral sequence degenerates at $E_2$, and we have $H^q(G^{\bullet} \times (G \times V), \u{\T}) = E^{0, q}_2$. Let $q$ be a positive integer. In this case, we have $E^{0, q}_1 \cong H^{q+1}(G, \Z)$. Under this isomorphism, we can see that $\phi_0^*\d_0^*c = 0$ for a class $c \in E^{0, q}_1$. We also have $\phi_0^*\d_1^*c = c$, because $\d_1 \circ \phi_0 = \id$. Thus, $E^{0, q}_2 = 0$ for $q > 0$. It is direct to see $E^{0, 0}_2 = H^0(G^{\bullet} \times (G \times V), \u{\T}) = C^{\infty}(V, \T)$, which completes the proof. 
\end{proof}

\begin{thm} \label{thm:quotient_simplicial_uT}
Let $G$ be a Lie group acting on a smooth manifold $M$. We assume that the action is free and locally trivial, and that the quotient space $M/G$ is a smooth manifold in such a way that the projection map $q : M \to M/G$ is smooth. For a non-negative integer $m$ the projection map induces an isomorphism of groups
\begin{eqnarray*}
q^* : H^m(M/G, \u{\T}) \to H^m(G^{\bullet} \times M, \u{\T}).
\end{eqnarray*}
\end{thm}

\begin{proof}
We use the spectral sequence in Lemma \ref{lem:Leray_type_ss_uT}. By the hypothesis, any point $\bar{x} \in M/G$ has a neighborhood $V$ such that $q^{-1}(V)$ is equivariantly isomorphic to $G \times V$. We can take $V$ to be a contractible open subset. By means of Lemma \ref{lem:cohomology_vanishing_uT}, the sheaf $\u{\mathcal{X}}^0$ is identified with $\u{\T}_{M/G}$. We also have $\u{\mathcal{X}}^q = 0$ for $q > 0$. Thus the degeneration of the spectral sequence at $E_2$ yields an isomorphism $H^m(M/G, \u{\T}_{M/G}) \cong H^m(G^{\bullet} \times M, \u{\T})$. We can see that this isomorphism is composed of $H^m(M/G, \u{\T}_{M/G}) \cong H^m(\{ e \}^{\bullet} \times M/G, \u{\T})$ given in Lemma \ref{lem:G_trivial} and $q^* : H^m(\{ e \}^{\bullet} \times M/G, \u{\T}) \to H^m(G^{\bullet} \times M, \u{\T})$ induced from the simplicial map $q : G^{\bullet} \times M \to \{ e \}^{\bullet} \times M/G$.
\end{proof}

\begin{rem}
As is mentioned in Remark \ref{rem_EC}, if $G$ is compact, then $H^m(G^{\bullet} \times M, \u{\T}) \cong H^{m+1}_G(M, \Z)$. Note that $H^m(M/G, \u{\T}) \cong H^{m+1}(M/G, \Z)$ provided that $M/G$ is a smooth manifold. As is well-known \cite{A-B}, if $G$ acts on $M$ freely, then $H^{m+1}_G(M, \Z) \cong H^{m+1}(M/G, \Z)$. Assembling these facts, we obtain an easier proof of Theorem \ref{thm:quotient_simplicial_uT} in the case that $G$ is compact and $m$ is positive.
\end{rem}

\subsection{The Deligne cohomology on the quotient space}

We denote by $A^q(M)_{cl}$ the group of closed $q$-forms on $M$.

\begin{lem} \label{lem:sEDC_to_forms}
Let $N$ be a positive integer, and $G$ a Lie group acting on a smooth manifold $M$. The group $H^N(G^{\bullet} \times M, \F(N))$ fits into the exact sequence
$$
\xymatrix@C=10pt@R=10pt{
0 \ar[r] &
H^N(G^{\bullet} \times M, \T) \ar[r] &
H^N(G^{\bullet} \times M, \F(N)) \ar[r]^-d & 
A^{N+1}(M)^G_{cl, bas} \ar[d] \\ 
& & &
H^{N+1}(G^{\bullet} \times M, \T),
} 
$$
where $H^m(G^{\bullet} \times M, \T)$ is the hypercohomology group of the constant simplicial sheaf $\T$, and $A^{N+1}(M)^G_{cl,bas} = Ker\{ \d : A^{N+1}(M)_{cl} \to A^{N+1}(G \times M)_{cl} \}$.
\end{lem}

\begin{proof}
We have the following short exact sequence of complexes of simplicial sheaves on $G^{\bullet} \times M$:
\begin{eqnarray}
\xymatrix@C=10pt{
0 \ar[r] &
\{ \u{\T} \to \u{A}^1 \to \cdots \to \u{A}^N_{cl} \} \ar[r] &
\F(N) \ar[r]^-{d} &
\{ 0 \to \cdots \to 0 \to \u{A}^{N+1}_{cl} \} \ar[r] & 0,
}
\label{exact_seq:Deligne_to_forms}
\end{eqnarray}
where $\u{A}^q_{cl}$ is the simplicial sheaf on $G^{\bullet} \times M$ given by the sheaf of germs of closed $q$-forms on each $G^i \times M$. The \Poincare lemma \cite{Bry1} induces a quasi-isomorphism
$$
\{ \T \to 0 \to \cdots \to 0 \} \to
\{ \u{\T} \to \u{A}^1 \to \cdots \to \u{A}^N_{cl} \}.
$$
The \Poincare lemma also allows us to use the complex of simplicial sheaves $(\u{A}^{*+N+1}, d)$ as a resolution of $\u{A}^{N+1}_{cl}$. Since each $G^i \times M$ is assumed to admit a partition of unity by our convention, we obtain
$$
H^m(G^{\bullet} \times M, 0 \to \cdots \to 0 \to \u{A}^{N+1}_{cl}) =
\left\{
\begin{array}{cl}
0, & (0 \le m < N), \\
A^{N+1}(M)^G_{cl,bas}, & (m = N). \\
\end{array}
\right.
$$
Now the long exact sequence associated with (\ref{exact_seq:Deligne_to_forms}) leads to the lemma.
\end{proof}

\begin{lem} \label{lem:quotient_simplicial_T}
Let $G$ be a Lie group acting on a smooth manifold $M$. We assume that the action is free and locally trivial. For a non-negative integer $m$ the projection map $q : M \to M/G$ induces an isomorphism of groups
\begin{eqnarray*}
q^* : H^m(M/G, \T) \to H^m(G^{\bullet} \times M, \T).
\end{eqnarray*}
\end{lem}

\begin{proof}
The argument here is the same as that in the proof of Theorem \ref{thm:quotient_simplicial_uT}. First, by the same method as in Lemma \ref{lem:Leray_type_ss_uT}, we have a spectral sequence converging to a graded quotient of the cohomology group $H^{p+q}(G^{\bullet} \times M, \T)$. Its $E_2$-terms are given by
$$
E^{p, q}_2 = H^p(M/G, \u{\mathcal{H}}^q),
$$
where $\u{\mathcal{H}}^q$ is the sheaf on $M/G$ associated with the presheaf given by the assignment of $H^q(G^{\bullet} \times q^{-1}(V), \T)$ to an open subset $V \subset M/G$. Second, by a method similar to that used in the proof of Lemma \ref{lem:cohomology_vanishing_uT}, we have 
$$
H^m(G^{\bullet} \times (G \times W), \T) =
\left\{
\begin{array}{cl}
\T, & \quad (m = 0), \\
0, & \quad (m > 0),
\end{array}
\right.
$$
where $W$ is a contractible set, and the Lie group $G$ acts on $G \times W$ by the left translation on $G$ and the trivial action on $W$. Finally, under the assumption in the current lemma, we identify $\u{\mathcal{H}}^0$ with the constant sheaf $\T$ on $M/G$ and $\u{\mathcal{H}}^q$ with the trivial sheaf $0$ for $q > 0$,  Then the degeneration of the spectral sequence completes the proof.
\end{proof}

\begin{thm} \label{thm:Deligne_coh_quotient}
Let $G$ and $M$ be as in Theorem \ref{thm:quotient_simplicial_uT}. For a non-negative integer $N$ the projection map $q : M \to M/G$ induces an isomorphism of groups
\begin{eqnarray*}
q^* : H^N(M/G, \F(N)) \to H^N(G^{\bullet} \times M, \F(N)).
\end{eqnarray*}
\end{thm}

\begin{proof}
Let $N$ be positive. The simplicial map $q : G^{\bullet} \times M \to \{ e \}^{\bullet} \times M/G$ induces a homomorphism between the exact sequences in Lemma \ref{lem:sEDC_to_forms}:
$$
\xymatrix{
0 \ar[d] &
0 \ar[d] \\
H^N(M/G, \T) \ar[r]^{q^*} \ar[d]&
H^N(G^{\bullet} \times M, \T) \ar[d] \\
H^N(M/G, \F(N)) \ar[r]^{q^*} \ar[d] &
H^N(G^{\bullet} \times M, \F(N)) \ar[d] \\
A^{N+1}(M/G)_{cl} \ar[r]^{q^*} \ar[d] &
A^{N+1}(M)^G_{cl,bas} \ar[d] \\
H^{N+1}(M/G, \T) \ar[r]^{q^*} &
H^{N+1}(G^{\bullet} \times M, \T),
}
$$
where the cohomology groups on $\{ e \}^{\bullet} \times (M/G)$ are identified with those on $M/G$ by Lemma \ref{lem:G_trivial}, and $A^{N+1}(M/G)^{\{ e \}}_{cl,bas}$ is identified with $A^{N+1}(M/G)_{cl}$. Note that $q : M \to M/G$ is a principal $G$-bundle by the assumption. Note also that $A^{N+1}(M)^G_{cl,bas}$ coincides with the group of closed basic $(N+1)$-forms on $M$ with respect to $q : M \to M/G$. Hence $q^* : A^{N+1}(M/G)_{cl} \to A^{N+1}(M)^G_{cl,bas}$ is an isomorphism. Now, Lemma \ref{lem:quotient_simplicial_T} and the five lemma establish the theorem in the case of $N$ positive. If $N = 0$, then the group $H^0(G^{\bullet} \times M, \F(0)) = H^0(G^{\bullet} \times M, \bF(0))$ is isomorphic to the group of $G$-invariant $\T$-valued smooth functions on $M$. Thus, the cohomology group is isomorphic to $H^0(M/G, \F(0)) = C^{\infty}(M/G, \T)$ under the assumption.
\end{proof}

Since $\bF(N)$ is obtained as the quotient $\bF(N) = \F(N) / F^1\!\F(N)$, we have a short exact sequence of complexes of simplicial sheaves:
\begin{eqnarray}
0 \longrightarrow 
F^1\!\F(N) \longrightarrow 
\F(N) \stackrel{\varphi}{\longrightarrow}
\bF(N) \longrightarrow 0.
\label{exact_seq:1st_filtration}
\end{eqnarray}
This induces a long exact sequence:
$$
\xymatrix@C=7pt@R=2pt{
&
\quad \quad \quad \quad \quad \quad \quad \quad \cdots \ar[r] &
H^{m-1}(G^{\bullet} \! \times \! M, \F(N)) \ar[r]^{\varphi} &
H^{m-1}(G^{\bullet} \! \times \! M, \bF(N)) \\
\ar[r]^-{\beta} &
H^m(G^{\bullet} \! \times \! M, F^1\!\F(N)) \ar[r] &
H^m(G^{\bullet} \! \times \! M, \F(N)) \ar[r]^{\varphi} &
H^m(G^{\bullet} \! \times \! M, \bF(N)) \\
\ar[r]^-{\beta} &
H^{m+1}(G^{\bullet} \! \times \! M, F^1\!\F(N)) \ar[r] &
H^{m+1}(G^{\bullet} \! \times \! M, \F(N)) \ar[r] &
\cdots. \quad \quad \quad  \quad \quad \quad
}
$$
We denote the Bockstein homomorphism by 
$$
\beta : 
H^m(G^{\bullet} \times M, \bF(N)) \longrightarrow
H^{m+1}(G^{\bullet} \times M, F^1\!\F(N)).
$$ 

Because $G^i \times M$ is assumed to admit a partition of unity for each $i$, the cohomology $H^m(G^{\bullet} \times M, F^1\!\F(N))$ is computed as the $m$th cohomology of the double complex $(L^{i, j}, \d, \til{d})$ given by
\begin{eqnarray}
L^{i, j} = 
\left\{
\begin{array}{cl}
F^1\!A^j(G^i \times M), & \quad (1 \le j \le N), \\
0, & \quad \mbox{otherwise}.
\end{array}
\right. \label{double_complex:F1FN}
\end{eqnarray}
It is clear that $H^m(G^{\bullet} \times M, F^1\!\F(N)) = 0$ for $m < N$. 

\begin{rem}
The image of $\beta : H^N(G^{\bullet} \times M, \bF(N)) \to H^{N+1}(G^{\bullet} \times M, F^1\!\F(N))$ plays the role of ``equivariant extensions'' in constructing a map from the equivariant smooth Deligne cohomology to the equivariant de Rham cohomology \cite{Go1}.
\end{rem}

\section{Equivariant circle bundles with connection}
\label{sec:ecbc}

In this section, we study relationships between equivariant circle bundles (with connection) and circle bundles (with connection) over the quotient space. Although the relationships can be seen directly, we make use of the results in the previous section.

We first consider the case of equivariant circle bundles \textit{without} connection.

\begin{prop} \label{prop:red_equiv_circle_bundle}
Let $G$ and $M$ be as in Theorem \ref{thm:quotient_simplicial_uT}. 

(a) For a $G$-equivariant principal $\T$-bundles $P$ over $M$, there exists a principal $\T$-bundle $\bar{P}$ over $M/G$ whose pull-back under the projection $q : M \to M/G$ is equivariantly isomorphic to $P$. 

(b) The isomorphism class of such $\bar{P}$ is unique.
\end{prop}

\begin{proof}
Recall Proposition \ref{prop:classify} (a) and Proposition \ref{prop:classify_equiv} (a). Let $c \in H^1(G^{\bullet} \times M, \u{\T})$ be the cohomology class corresponding to the equivariant isomorphism class of $P$. Since $q^* : H^1(M/G, \u{\T}) \to H^1(G^{\bullet} \times M, \u{\T})$ is an isomorphism by Theorem \ref{thm:quotient_simplicial_uT}, we put $\bar{c} = (q^*)^{-1}(c)$. Let $\bar{P}$ be a principal $\T$-bundle over $M/G$ which is classified by $\bar{c} \in H^1(M/G, \u{\T})$. Because $q^*(\bar{c}) = c$, the pull-back of $\bar{P}$ under the projection $q : M \to M/G$ is equivariantly isomorphic to $P$. 
\end{proof}

Note that we can directly construct a principal $\T$-bundle $\bar{P}$ such that  $q^*\bar{P}$ is equivariantly isomorphic to $P$ as follows. If the action of $G$ on $M$ is as in Theorem \ref{thm:quotient_simplicial_uT}, then so is the action of $G$ on $P$. Hence the quotient space $P/G$ gives rise to a principal $\T$-bundle over $M/G$. Clearly, the pull-back of $P/G \to M/G$ under $q : M \to M/G$ is equivariantly isomorphic to $P \to M$.

\smallskip

Next we consider the case of equivariant circle bundles \textit{with} connection. Recall that the cohomology $H^m(G^{\bullet} \times M, F^1\!\F(1))$ can be given by
$$
H^m(G^{\bullet} \times M, F^1\!\F(1)) =
\frac{Ker\{ \d : F^1\!A^1(G^{m-1} \times M) \to F^1\!A^1(G^m \times M)\}}
{Im\{ \d : F^1\!A^1(G^{m-2} \times M) \to F^1\!A^1(G^{m-1} \times M)\}}.
$$

\begin{lem} \label{vanising_H1F1F1}
$H^1(G^{\bullet} \times M, F^1\!\F(1)) = 0$.
\end{lem}

\begin{proof}
Because $F^1\!A^1(M) = 0$, this lemma is clear.
\end{proof}

\begin{prop} \label{prop:reduction_T_bundle_connection}
Let $G$ and $M$ be as in Theorem \ref{thm:quotient_simplicial_uT}, $(P, \theta)$ a $G$-equivariant principal $\T$-bundle with $G$-invariant connection over $M$, and $c \in H^1(G^{\bullet} \times M, \bF(1))$ the class corresponding to the equivariant isomorphism class of $(P, \theta)$. 

(a) There exists a principal $\T$-bundle with connection $(\bar{P}, \bar{\theta})$ over $M/G$ whose pull-back under the projection $q : M \to M/G$ is equivariantly isomorphic to $(P, \theta)$ if and only if $\beta(c) = 0$ in $H^2(G^{\bullet} \times M, F^1\!\F(1))$.

(b) The isomorphism class of such $(\bar{P}, \bar{\theta})$ is unique.
\end{prop}

\begin{proof}
By Theorem \ref{thm:Deligne_coh_quotient} and Lemma \ref{vanising_H1F1F1}, the short exact sequence (\ref{exact_seq:1st_filtration}) gives the following exact sequence
$$
0 \longrightarrow
H^1(M/G, \F(1)) \stackrel{\varphi \circ q^*}{\longrightarrow} 
H^1(G^{\bullet} \times M, \bF(1)) \stackrel{\beta}{\longrightarrow} 
H^2(G^{\bullet} \times M, F^1\!\F(1)).
$$
There exists a cohomology class $\bar{c} \in H^1(M/G, \F(1))$ such that $\varphi \circ q^*(\bar{c}) = c$ if and only if $\beta(c) = 0$. Let $(\bar{P}, \bar{\theta})$ be a principal $\T$-bundle with connection over $M/G$ classified by $\bar{c}$. Since $\varphi \circ q^*(\bar{c}) = c$, the pull-back of $(\bar{P}, \bar{\theta})$ under $q : M \to M/G$ is equivariantly isomorphic to $(P, \theta)$. Because $\varphi \circ q^*$ is injective, the isomorphism class of $(\bar{P}, \bar{\theta})$ is unique.
\end{proof}

Let $\g$ be the Lie algebra of $G$, and $\langle \ | \ \rangle : \g \otimes \g^* \to \R$ the natural contraction. By the (co)adjoint, the Lie group $G$ acts on an element $f \in \g^*$ by $\langle X | \Ad_{g} f \rangle = \langle \Ad_{g^{-1}} X| f \rangle$. 

\begin{lem} \label{lem:iso_H2F1F1}
There exists an isomorphism
\begin{eqnarray}
H^2(G^{\bullet} \times M, F^1\!\F(1)) \cong
\left\{ f : M \to \g^* |\ 
g^*f = \Ad_g f \ \mbox{for all} \ g \in G \right\}.
\label{expression_cohomology_F1F1}
\end{eqnarray}
\end{lem}

\begin{proof}
Because $F^1\!A^1(M) = 0$, we have
$$
H^2(G^{\bullet} \times M, F^1\!\F(1)) =
\{ \alpha \in  F^1\!A^1(G \times M) |\ \d \alpha = 0 \}.
$$
Let $\alpha$ be an element in $F^1\!A^1(G \times M)$. Note that for any tangent vector $V \in T_xM$ we have $\alpha((g, x); V) = 0$. By a computation, we can see that the cocycle condition $\d \alpha = 0$ is equivalent to the following conditions:
\begin{eqnarray*}
\alpha((g_2, x); g_2X) & = & \alpha((g_1g_2, x); g_1g_2X), \\
\alpha((g_1, g_2x); g_1X) & = & \alpha((g_1, g_2x); g_1Xg_2),
\end{eqnarray*}
where a tangent vector at $g \in G$ is expressed as $gX \in T_gG$ by an element $X \in T_eG = \g$. Thus, the isomorphism (\ref{expression_cohomology_F1F1}) is induced by the assignment to $\alpha$ of the map $f : M \to \g^*$ defined by $\langle X | f(x) \rangle = \alpha((e, x); X)$. The inverse homomorphism is given by the assignment to $f : M \to \g^*$ of the 1-form $\alpha$ defined by $\alpha((g, x); gX \oplus V) = \langle X | f(x) \rangle$.
\end{proof}

For a $G$-invariant connection $\theta$ on a $G$-equivariant principal $\T$-bundle $P$ over $M$, we define a map $\mu : M \to \g^*$ by
\begin{eqnarray}
\langle X | \mu(x) \rangle = 
\frac{-1}{2\pi\i}\theta(p; X^*), \label{moment}
\end{eqnarray}
where $p \in P$ is a point lying on the fiber of $x$, and $X^* \in T_pP$ is the tangent vector generated by the infinitesimal action of $X \in \g$. Since the $G$-action on $P$ commutes with the right $\T$-action on $P$, the map $\mu$ is well-defined. We call the map $\mu$ the \textit{moment} \cite{B-V} associated with $(P, \theta)$.

\begin{lem}[\cite{Go1}] \label{lem:beta_c:degree1}
Let $(P, \theta)$ be a $G$-equivariant $\T$-bundle with $G$-invariant connection over $M$, and $c \in H^1(G^{\bullet} \times M, \bF(1))$ the cohomology class that classifies $(P, \theta)$. Under the isomorphism (\ref{expression_cohomology_F1F1}), the image $\beta(c) \in H^2(G^{\bullet} \times M, F^1\!\F(1))$ is identified with the moment $\mu : M \to \g^*$ associated with $(P, \theta)$. 
\end{lem}

By the help of Lemma \ref{lem:beta_c:degree1}, we can directly construct a principal $\T$-bundle $\bar{P}$ with a connection $\bar{\theta}$ such that $q^*(\bar{P}, \bar{\theta})$ is equivariantly isomorphic to $P$. As we see, $P$ gives a principal $\T$-bundle $\bar{P} \to M/G$ by setting $\bar{P} = P/G$. By Lemma \ref{lem:beta_c:degree1}, the condition $\beta(c) = 0$ in Proposition \ref{prop:reduction_T_bundle_connection} is equivalent to the vanishing of the moment $\mu$ associated with $(P, \theta)$. It is clear by (\ref{moment}) that the condition $\mu \equiv 0$ is the necessary and sufficient condition for the $G$-invariant connection $\theta$ on $P$ to descend to induce a connection $\bar{\theta}$ on the principal $\T$-bundle $P/G \to M/G$. The pull-back of $(\bar{P}, \bar{\theta})$ under $q : M \to M/G$ is equivariantly isomorphic to $(P, \theta)$.

\section{Equivariant gerbe with connection}
\label{sec:egc}

\begin{thm} \label{thm:red_equiv_gerbe}
Let $G$ and $M$ be as in Theorem \ref{thm:quotient_simplicial_uT}. 

(a) For a $G$-equivariant gerbe $\sC$ over $M$, there exists a gerbe $\bsC$ over $M/G$ whose pull-back under the projection $q : M \to M/G$ is equivariantly isomorphic to $\sC$. 

(b) The isomorphism class of such $\bsC$ is unique.
\end{thm}

\begin{proof}
The proof is the same as that of Proposition \ref{prop:red_equiv_circle_bundle}. By Proposition \ref{prop:classify_equiv} (b), we have a class $c \in H^2(G^{\bullet} \times M, \u{\T})$ that classifies the equivariant isomorphism class of $\sC$. Since $q^* : H^2(M/G, \u{\T}) \to H^2(G^{\bullet} \times M, \u{\T})$ is an isomorphism by Theorem \ref{thm:quotient_simplicial_uT}, we put $\bar{c} = (q^*)^{-1}(c)$. Let $\bsC$ be a gerbe over $M/G$ corresponding to $\bar{c} \in H^2(M/G, \u{\T})$ by Proposition \ref{prop:classify} (b). Because $q^*(\bar{c}) = c$, the pull-back of $\bsC$ under the projection $q : M \to M/G$ is equivariantly isomorphic to $\sC$. 
\end{proof}

Since a gerbe does not have a ``total space'', the construction of such a gerbe $\bsC$ as in Theorem \ref{thm:red_equiv_gerbe} would be not as direct as in the case of circle bundles. For \textit{equivariant bundle gerbes}, which are closely related to equivariant gerbes, we have some constructions of bundle gerbes over the quotient space \cite{Gaw-R,Ma-S,Me}.

\begin{thm} \label{thm:reduction_gerbe_connection}
Let $G$ and $M$ be as in Theorem \ref{thm:quotient_simplicial_uT}, $(\sC, \Co, K)$ a $G$-equivariant gerbe with $G$-invariant connective structure and $G$-invariant curving over $M$, and $c \in H^2(G^{\bullet} \times M, \bF(2))$ the class corresponding to the equivariant isomorphism class of $(\sC, \Co, K)$. 

(a) There exists a gerbe with connective structure and curving $(\bsC, \bCo, \bK)$ over $M/G$ whose pull-back under the projection $q : M \to M/G$ is equivariantly isomorphic to $(\sC, \Co, K)$ if and only if $\beta(c) = 0$ in $H^3(G^{\bullet} \times M, F^1\!\F(2))$.

(b) The isomorphism classes of such $(\bsC, \bCo, \bK)$ are in one to one correspondence with $Coker\{ \beta : H^1(G^{\bullet} \times M, \bF(2)) \to H^2(G^{\bullet} \times M, F^1\!\F(2)) \}$.
\end{thm}

\begin{proof}
By Theorem \ref{thm:Deligne_coh_quotient} and (\ref{exact_seq:1st_filtration}), we have the following exact sequence:
$$
\xymatrix@R=2pt{
&
H^1(G^{\bullet} \! \times \! M, \bF(2)) \ar[r]^-{\beta} &
H^2(G^{\bullet} \! \times \! M, F^1\!\F(2)) \ar[r] &
H^2(M/G, \F(2)) \\
\ar[r]^-{\varphi \circ q^*} &
H^2(G^{\bullet} \! \times \! M, \bF(2)) \ar[r]^-{\beta} &
H^3(G^{\bullet} \! \times \! M, F^1\!\F(2)). &
}
$$
There is a class $\bar{c} \in H^2(M/G, \F(2))$ such that $\varphi \circ q^*(\bar{c}) = c$ if and only if $\beta(c) = 0$. By Proposition \ref{prop:classify_connection}, there is a gerbe with connective structure and curving $(\bsC, \bCo, \bK)$ over $M/G$ classified by $\bar{c}$. Since $\varphi \circ q^*(\bar{c}) = c$, the pull-back of $(\bsC, \bCo, \bK)$ under $q$ is equivariantly isomorphic to $(\sC, \Co, K)$. Hence (a) is proved. By the exact sequence above, we have a bijection between the set of $\bar{c} \in H^2(M/G, \F(2))$ such that $\varphi \circ q^*(\bar{c}) = c$ and the cokernel of $\beta : H^1(G^{\bullet} \times M, \bF(2)) \to H^2(G^{\bullet} \times M, F^1\!\F(2))$, which leads to (b). 
\end{proof}

We notice that, although we may have distinct isomorphism classes of such $(\bsC, \bCo, \bK)$ as in Theorem \ref{thm:reduction_gerbe_connection}, the isomorphism class of $\bsC$ is unique. This is a consequence of Theorem \ref{thm:red_equiv_gerbe}. We also notice that we can characterize the 3-curvature (\cite{Bry1}) of $(\bsC, \bCo, \bK)$ by the unique 3-form $\bar{\Omega} \in A^3(M/G)$ such that $q^*\bar{\Omega} = \Omega$, where $\Omega \in A^3(M)$ is the 3-curvature of $(\sC, \Co, K)$ over $M$.

\begin{lem} \label{lem:iso_H3F1F2}
There is an isomorphism $H^3(G^{\bullet} \times M, F^1\!\F(2)) \cong \mathcal{Z}/\mathcal{B}$. Here $\mathcal{Z}$ and $\mathcal{B}$ are defined by
\begin{eqnarray}
\mathcal{C} & = & 
A^1(M, \g^*) \oplus A^0(G \times M, \g^*), \\
\mathcal{Z} & = &
\left\{
(E, \zeta) \in \mathcal{C} \Big|
\begin{array}{l}
E(gx; gV) - Ad_gE(x; V) = d_M\zeta((g, x); V), \\
\Ad_g\zeta(h, x) - \zeta(gh, x) + \zeta(g, hx) = 0. \\
\end{array}
\right\}, \\ 
\mathcal{B} & = &
\left\{
(E, \zeta) \in \mathcal{C} \Bigg|
\begin{array}{l}
f \in A^0(M, \g^*), \\
E = d f, \\
\zeta(g, x) = f(gx) - Ad_g f(x).
\end{array}
\right\}, \label{subspace:coboundary}
\end{eqnarray}
where $d_M$ is the exterior differential in the direction of $M$.
\end{lem}

\begin{proof}
Recall that $H^3(G^{\bullet} \times M, F^1\!\F(2))$ is the cohomology of (\ref{double_complex:F1FN}). Let $[\alpha, \beta]$ be a class in the cohomology, where $\alpha \in F^1\!A^2(G \times M)$ and $\beta \in F^1\!A^1(G^2 \times M)$. We define $H\!\beta \in  F^1\!A^1(G \times M)$ by setting $H\!\beta((g, x); gX \oplus V) = \beta((g, e, x); 0 \oplus X \oplus 0)$. Note that cocycles $(\alpha, \beta)$ and $(\alpha - d H\!\beta, \beta + \d H\!\beta)$ induce the same cohomology class. Now we define $E \in A^1(M, \g^*)$ and $\zeta \in A^0(G \times M, \g^*)$ by
\begin{eqnarray*}
\langle X | E(x; V) \rangle
& = &
(\alpha - d H\!\beta)((e, x); X \oplus 0, 0 \oplus V), \\
\langle X | \zeta(g, x) \rangle
& = &
(\beta + \d H\!\beta)((e, g, x); X \oplus 0 \oplus 0).
\end{eqnarray*}
By a lengthy calculation, we can see that the cocycle condition for $(\alpha, \beta)$ implies that $(E, \zeta)$ belongs to $\mathcal{Z}$. We here consider the case that $(\alpha, \beta) = (- d\gamma, \d \gamma)$ for an element $\gamma \in F^1\!A^1(G \times M)$. In this case, we define $f \in A^0(M, \g^*)$ by
$$
\langle X | f(x) \rangle = \gamma((e, x); X \oplus 0).
$$
By using $f$ above, we can express the $(E, \zeta) \in \mathcal{Z}$ defined by $(\alpha, \beta)$ as in (\ref{subspace:coboundary}). Therefore we obtain a homomorphism $\Phi : H^3(G^{\bullet} \times M, F^1\!\F(2)) \to \mathcal{Z}/\mathcal{B}$ by setting $\Phi([\alpha, \beta]) = [E, \zeta]$. To show that $\Phi$ is an isomorphism, we give the inverse homomorphism. For $(E, \zeta) \in \mathcal{C}$ we define $\alpha \in F^1\!A^2(G \times M)$ and $\beta \in F^1\!A^1(G^2 \times M)$ by
\begin{eqnarray*}
\alpha((g, x); gX \oplus V, gX' \oplus V')
& = &
\langle X | E(x; V') \rangle - 
\langle X' | E(x; V) \rangle \\
& + &
\langle [X, X'] | \zeta(e, x) \rangle, \\
\beta((g_1, g_2, x); g_1X_1 \oplus g_2 X_2 \oplus V)
& = &
\langle X_1 | \zeta(g_2, x) \rangle.
\end{eqnarray*}
If $(E, \zeta)$ belongs to $\mathcal{Z}$, then $(\alpha, \beta)$ is a cocycle. We here suppose that $(E, \zeta)$ can be expressed as in (\ref{subspace:coboundary}) by a function $f \in A^0(M, \g^*)$. In this case, we define $\gamma \in F^1\!A^1(G \times M)$ by
$$
\gamma((g, x); gX \oplus V) = \langle X | f(x) \rangle.
$$
We can verify that the $(\alpha, \beta)$ defined by $(E, \zeta)$ is expressed as $(-d \gamma, \d \gamma)$. Therefore we obtain a homomorphism $\Psi : \mathcal{Z}/\mathcal{B} \to H^3(G^{\bullet} \times M, F^1\!\F(2))$ by setting $\Psi([E, \zeta]) = [\alpha, \beta]$. Note that, if $(E, \zeta) \in \mathcal{Z}$, then we have 
$$
\langle [X, X'] | \zeta(e, x) \rangle 
=
\langle X | E(x; {X'}^* ) \rangle -
\langle X | d\zeta((e, x); X' \oplus 0) \rangle.
$$
Thus, we can see that $\Psi$ is the inverse of $\Phi$.
\end{proof}

An equivariant gerbe with connective structure and curving does not determine a representative $(E, \zeta) \in \mathcal{Z}$ of $\beta(c)$, in general. So we do not spell out here the general formula of the representative. (In \cite{Bry2}, a computation of the term $E$ can be found.)

\begin{lem} 
There is an isomorphism
$$
H^2(G^{\bullet} \times M, F^1\!\F(2)) \cong
\{ f \in A^0(M, \g^*) |\ 
df = 0, \ g^*f = \Ad_g f \ \mbox{for all} \ g \in G \}.
$$
In particular, if $M$ is connected, then there is an isomorphism
$$
H^2(G^{\bullet} \times M, F^1\!\F(2)) \cong
\{ f \in \g^* |\ f = \Ad_g f \ \mbox{for all} \ g \in G \}.
$$
\end{lem}

\begin{proof}
Because $F^1\!A^1(M) = F^1\!A^2(M) = 0$, we have
$$
H^2(G^{\bullet} \times M, F^1\!\F(2)) =
\left\{
\alpha \in F^1\!A^1(G \times M) |\ \d \alpha = 0, \ d \alpha = 0
\right\}.
$$
By Lemma \ref{lem:iso_H2F1F1}, an element $\alpha \in F^1\!A^1(G \times M)$ such that $\d \alpha = 0$ is uniquely expressed as $\alpha = \langle g^{-1}dg | f \rangle$, where $g^{-1}dg$ is the left invariant Maurer-Cartan form on $G$ and $f \in A^0(M, \g)$ satisfies $g^*f = \Ad_gf$. The differential of $\alpha$ is 
$$
d \alpha = 
\langle dg^{-1} \wedge dg | f \rangle - 
\langle g^{-1}dg | df \rangle \in F^1\!A^2(G \times M). 
$$
Note that $\langle dg^{-1} \wedge dg | f \rangle \in F^2\!A^2(G \times M)$, while $\langle g^{-1}dg | df \rangle \not\in F^2\!A^2(G \times M)$. So the condition $d\alpha = 0$ is equivalent to $df = 0$ and $\langle [X, Y] | f \rangle = 0$ for all $X, Y \in \g$. The last condition follows from $d f = 0$ and $g^*f = \Ad_g f$.
\end{proof}

The proof above implies that, if $\g = [\g, \g]$, then $H^2(G^{\bullet} \times M, F^1\!\F(2))= 0$.

\begin{cor}
If the Lie algebra $\g$ of $G$ is such that $[\g, \g] = \g$, then the isomorphism class of $(\bsC, \bCo, \bK)$ in Theorem \ref{thm:reduction_gerbe_connection} (a) is unique.
\end{cor}

In contrast with the above, if $G$ contains tori as a center, then the group $H^2(G^{\bullet} \times M, F^1\!\F(2))$ is non-trivial. For example, let $T$ be a maximal torus of $SU(2)$. If we put $G = T$ and $M = SU(2)$, then we have $M/G \cong S^2$. In this case, we obtain $H^2(G^{\bullet} \times M, F^1\!\F(2)) \cong \R$. We also obtain by computations $H^1(G^{\bullet} \times M, \bF(2)) \cong \Z$ and $H^1(G^{\bullet} \times M, \F(2)) = 0$. As a result, we have $Coker\{ \beta : H^1(G^{\bullet} \times M, \bF(2)) \to H^2(G^{\bullet} \times M, F^1\!\F(2)) \} \cong \R/\Z$. 

\medskip

In general, the cohomology group $H^1(G^{\bullet} \times M, \bF(2))$ classifies the isomorphism classes of $G$-equivariant $\T$-bundle over $M$ with \textit{flat} connection \cite{Go1}. The image of $\beta : H^1(G^{\bullet} \times M, \bF(2)) \to H^2(G^{\bullet} \times M, F^1\!\F(2))$ consists of the moment maps associated with $G$-equivariant $\T$-bundles with flat connection over $M$.


\bigskip

\begin{acknowledgments}
I would like to thank T. Kohno for useful discussions and valuable suggestions, and M. Furuta for good advice. I would also like to thank Y. Hashimoto for reading earlier drafts and for helpful comments.
\end{acknowledgments}



\begin{flushleft}
Graduate school of Mathematical Sciences, University of Tokyo, \\
Komaba 3-8-1, Meguro-Ku, Tokyo, 153-8914 Japan. \\
e-mail: kgomi@ms.u-tokyo.ac.jp
\end{flushleft}

\end{document}